\documentclass[twocolumn, 10pt]{IEEEtran}

\usepackage{times}
\usepackage[latin1]{inputenc}
\usepackage[T1]{fontenc}
\usepackage[english]{babel}
\usepackage{balance}

\usepackage{amsmath, amsfonts, amssymb, array, vector, booktabs}
\usepackage{subeqnarray}
\usepackage{algorithmic, algorithm}
\usepackage{amstext}

\usepackage{tikz}

\usepackage{graphicx}
\usepackage{caption}
\usepackage{subcaption}
\usepackage{psfrag}
\usepackage{epsf,epsfig}

\usepackage{multirow}
\usepackage{url}
\usepackage{color}
\usepackage[normalem]{ulem}

\newtheorem{thm}{Theorem}[section]
\newtheorem{prop}[thm]{Proposition}
\newcommand{\maxi}{\mathop{\displaystyle \mbox{maximize}}}

\newcommand{\st}{\mathop{\displaystyle \mbox{subject to}}}

\begin{document}
%

\title{A Reinforcement Learning Approach to Power Control and Rate Adaptation in Cellular Networks}


\author{\IEEEauthorblockN{Euhanna Ghadimi, Francesco Davide Calabrese, Gunnar Peters, Pablo Soldati }
	\\
\IEEEauthorblockA{
Huawei Technologies Sweden AB, R\&D Center, Kista, Sweden\\
e-mail: \{firstname.lastname\}@huawei.com
}
}
\maketitle

\begin{abstract}
Optimizing radio transmission power and user data rates in wireless systems via power control requires an accurate and instantaneous knowledge of the system model. While this problem has been extensively studied in the literature, an efficient solution approaching optimality with the limited information available in practical systems is still lacking.
This paper presents a reinforcement learning framework for power control and rate adaptation in the downlink of a radio access network that closes this gap.
We present a comprehensive design of the learning framework that includes the characterization of the system state, the design of a general reward function, and the method to learn the control policy. System level simulations show that our design can quickly learn a power control policy that brings significant energy savings and fairness across users in the system.
\end{abstract}


%
\IEEEpeerreviewmaketitle
\section{Introduction}\label{sec:1}

Radio interference can drastically degrade the performance of radio systems when not properly dealt with. Interference mitigation has thereby played a major role in radio access networks and shall continue to do so as we move toward the $5^{\mathrm{th}}$ generation (5G) of mobile broadband systems. Compared to the 4G systems, however, interference in 5G networks is expected to have different behavior and characteristics due to an increasingly heterogeneous, multi-RAT and denser environment to the extent that conventional inter-cell interference coordination (ICIC) will become inadequate~\cite{Andrews:14}.

Interference mitigation can broadly be posed as a power control optimization problem which, under certain conditions, admits an optimal solution. In particular, distributed power control methods based on linear iterations to meet signal-to-interference-and-noise ratio (SINR) targets were proposed in~\cite{Zander:92, FoM:93}. An axiomatic framework for studying general power control iterations was proposed in~\cite{Yates:95} based on the so called \emph{standard interference functions}, and extended in~\cite{SoL:05, BoS:10}. Despite the elegant and insightful results available in the power control literature, their direct application to practical systems has been impaired by the dependency on some simplifying assumptions, such as the knowledge of the instantaneous channel gains to/from all user devices, full buffer traffic, etc.

Recent advances in the field of machine learning vastly expanded the class of problems that can be tackled with such methods and made new techniques available that can potentially change the way radio resource management (RRM) problems are solved in wireless systems. Within this field, Reinforcement Learning (RL) is arguably the most appropriate branch in which one or more agents solve a complex control problem by interacting with an environment, issuing control actions, and utilizing the feedbacks obtained from the environment to find optimal control policies  (see e.g.~\cite{SuB:98}).

Only recently, learning algorithms made their appearance in the context of wireless networking. Cognitive radios presented a natural ground for the application of learning methods, for instance for efficient spectrum sensing and shaping see, e.g.~\cite{Macaluso:13a, Macaluso:14}, for mitigating interference generated by multiple cognitive radios at the receivers of primary users, see e.g.~\cite{GaG:10,SCG:11}, and so on. A survey of learning methods suitable for cognitive radio networks can be found in~\cite{Gavrilovska:13} and references therein.
Distributed RL was proposed in~\cite{GaG:10}~to enable radio cells (i.e., the agents) to learn an efficient policy for controlling the aggregated interference generated by multiple neighboring cells on primary licensed users. The proposed algorithm exploits a Q-learning method based on a representation with tables and neural networks. Table-based Q-learning was also considered in~\cite{SCG:11} for interference mitigation in LTE heterogeneous radio cellular networks (HetNets) and in~\cite{SBG:15} for frequency- and time-domain inter-cell interference coordination in HetNets.

In this paper we apply RL for distributed downlink inter-cell power control and rate adaptation in the downlink of a radio access network. We begin by modeling the problem as a network utility maximization and solve it to optimality via Lagrange duality theory. While technically sound, several practical limitations refrain us from considering this optimization framework for applications to real systems.
To overcome the impracticalities of such approach, we propose an advanced RL framework which is particularly data-efficient. By carefully designing the features required to characterize the system state for multi-cell downlink power control and rate adaptation as well as the reward function to drive the behavior of the agent, the framework is able to quickly produce a policy for power control that brings energy saving as well as fairness across the users in the system.
The validity of this approach is verified in a fully LTE-A compliant event-driven system level simulator.
The results also demonstrate the flexibility of this approach which enables us to promote entirely different behaviors in terms of fairness and system performance by changing the parameters of the reward function.

The rest of the paper is organized as follows: Sec.~\ref{sec:2} introduces a system model for the rate and power control optimization problem solved in Sec.~\ref{sec:3}. Sec.~\ref{sec:4} introduces a RL architecture for downlink power control. Finally, Sec.~\ref{sec:numerical} and~\ref{sec:conclusion} present simulation results and final remarks, respectively.
\section{System model}\label{sec:2}
We consider a radio cellular system with $C$ cells, labelled $c=1,\dots, C$, each serving a set of users $\mathcal{N}_c$. Users in the system are labelled by $n=1,\dots,N$ where $N = \sum_c N_c$ and $N_c=\mid\mathcal{N}_c\mid$. We assume a frequency reuse-1 scheme where all cells operate within the same frequency bandwidth $W$ with maximum downlink transmission power $P_c^{\max}$, respectively. We assume the system bandwidth divided into $K$ equally sized time-frequency resource blocks. Within each cell, users are orthogonally scheduled in frequency bandwidth $W$ so that only inter-cell interference is considered.

Let $P_c$ denote the downlink transmission power budget used by cell $c$ in a given transmission time interval (TTI), $P_c^{\max}$ be the corresponding maximum power budget (i.e. $P_c\leq P_c^{\max}$), and  $\mathbf{p} = [P_1,\dots, P_C]^T$ be the network-wide vector of transmission power budgets for all cells in a TTI. 
Assuming each cell uniformly distributes it transmission power budget $P_c$ over the time-frequency resource blocks scheduled to users in the cell, each user within a cell is served with the same power level, i.e. $P_c/K$ for all users $n\in\mathcal{N}_c$.
%

Furthermore, we define $G_{n,c}$ as the channel gain between user $n$ and cell $c$ which takes into account for large-scale fading effects (i.e., pathloss and shadowing). Therefore, the average SINR experienced by user $n$ can be modelled as
\begin{align}\label{eqn:SINR}
\gamma_n(\mathbf{p}) = \frac{G_{n,c(n)}P_{c(n)}}{\sigma^2 + \sum_{c\neq c(n)} G_{n,c} P_c} \quad n = 1,\dots,N.
\end{align}
where $c(n)$ denotes the serving cell of user $n$. We would like to stress that this model is applied \emph{exclusively} to the optimization framework developed in Section~\ref{sec:3} but it is not explicitly needed for the learning framework described in Section~\ref{sec:4}.


\section{Optimal rate and power control} \label{sec:3}

We consider the problem of joint user rate and power control optimization in a multi-celluar radio network. We pose the problem as a network-wide utility maximization
\begin{align}
\hspace{-3mm}\begin{array}{lll}
\maxi & \sum_n  u_n(r_n) &\\
\st   & r_n \leq W_n\log_2(1+\gamma_n(\mathbf{p})) & n = 1,\dots, N,\\
      & \mathbf{p}\preceq \mathbf{p}^{\max} &
\end{array}\label{eqn:NUMprimal}
\end{align}
wherein $W_n = W/N_{c(n)}$ denotes the average amount of bandwidth scheduled for user $n$ assuming an equal share of frequency resources, and $r_n$ models the theoretically achievable user data rate according to the Shannon bound. Associated with each user $n$ is a utility function $u_n(\cdot)$, which describes the utility of the user to communicate at rate $r_n$ (cf.~\cite{Kelly:98}). We assume that $u_n$ is increasing and strictly concave, with $u_n\rightarrow -\infty$ as $r_n \rightarrow 0^+$. Therefore, problem~\eqref{eqn:NUMprimal} aims at optimizing the downlink transmission power of each cell so as to maximize a network-wide utility of the users in the system.

\subsection{Convexification  and optimal solution}
Problem~\eqref{eqn:NUMprimal} is not convex in the variables $\mathbf{p}$ due to inter-cell interference in the rate expression. However, with a suitable log-transformation of both constraints and variables one can obtain an equivalent convex formulation, cf.~\cite{PSE:08}. In particular, we define $\tilde{r}_n = \log(r_n)$ and $\tilde{P}_n = \log(P_n)$ and rewrite~\eqref{eqn:NUMprimal} as
\begin{align}
\hspace{-3mm}\begin{array}{ll}
\maxi & \sum_n  u_n(e^{\tilde{r}_n}) \\
\st   & \tilde{r}_n \leq \log(W_n)+ \log(\log_2(1+\gamma_n(e^{\tilde{\mathbf{p}}})))\; \forall n,\\
      & \tilde{\mathbf{p}}\preceq \log(\mathbf{p}^{\max}),
\end{array}\label{eqn:NUMtransf}
\end{align}
wherein constraints are obtained with a log-transformation and variable change applied to both sides of the inequalities in~\eqref{eqn:NUMprimal}.

\begin{prop}
The transformed problem~\eqref{eqn:NUMtransf} is jointly convex in $\tilde{\mathbf{r}}$ and $\tilde{\mathbf{p}}$. The optimal solution to~\eqref{eqn:NUMprimal} coincides with the one of~\eqref{eqn:NUMtransf}.
\end{prop}
\begin{IEEEproof}
The transformed capacity constraints are jointly convex in $\tilde{\mathbf{p}}$ and linear in $\tilde{\mathbf{r}}$, while the remaining power budget constraints are linear in $\tilde{\mathbf{p}}$. Therefore, the rest of the results follows from~\cite[Theorem~2]{PSE:08} and ~\cite[Corollary~1]{PSE:08}.
\end{IEEEproof}
An optimal and distributed solution to the rate and power control problem~\eqref{eqn:NUMprimal} can be found by solving~\eqref{eqn:NUMtransf} with standard Lagrange duality theory\footnote{The details are omitted due to space limitations but we refer to~\cite{PSE:08} for a similar approach.}. Examples of utility functions satisfying these conditions include the $\alpha$-fair utility functions, such as sum-log utility.

\subsection{Practical considerations}
While theoretically sound, solving problem~\eqref{eqn:NUMprimal} to optimality with the signaling and time limitations of conventional radio access networks would be infeasible. Firstly, solving problem~\eqref{eqn:NUMtransf} to optimality requires cross-cell knowledge of channel gains $G_{n,c}$ (or similar message passing see e.g.~\cite{PSE:08}). While this would create large signaling overhead, the channel aging  would render the information outdated before it serves its purpose. On the other hand, the user devices can in practice estimate only the channel gain $G_{n,c}$ for a small subset of interfering cells. Secondly, any change in the users and traffic distributions (e.g., user positions, number of users, number of resources used, etc.) would also need to be exchanged between the radio access nodes as it affects the achievable user data rate. Such changes are difficult to track, measure and expensive to be communicated in the system. These impairments have refrained the application of inter-cell interfere coordination via optimal power control in state of the art wireless systems.


\section{RL for Cellular Networks} \label{sec:4}

\begin{figure}[t!]
\centering{
\includegraphics[width=0.6\linewidth]{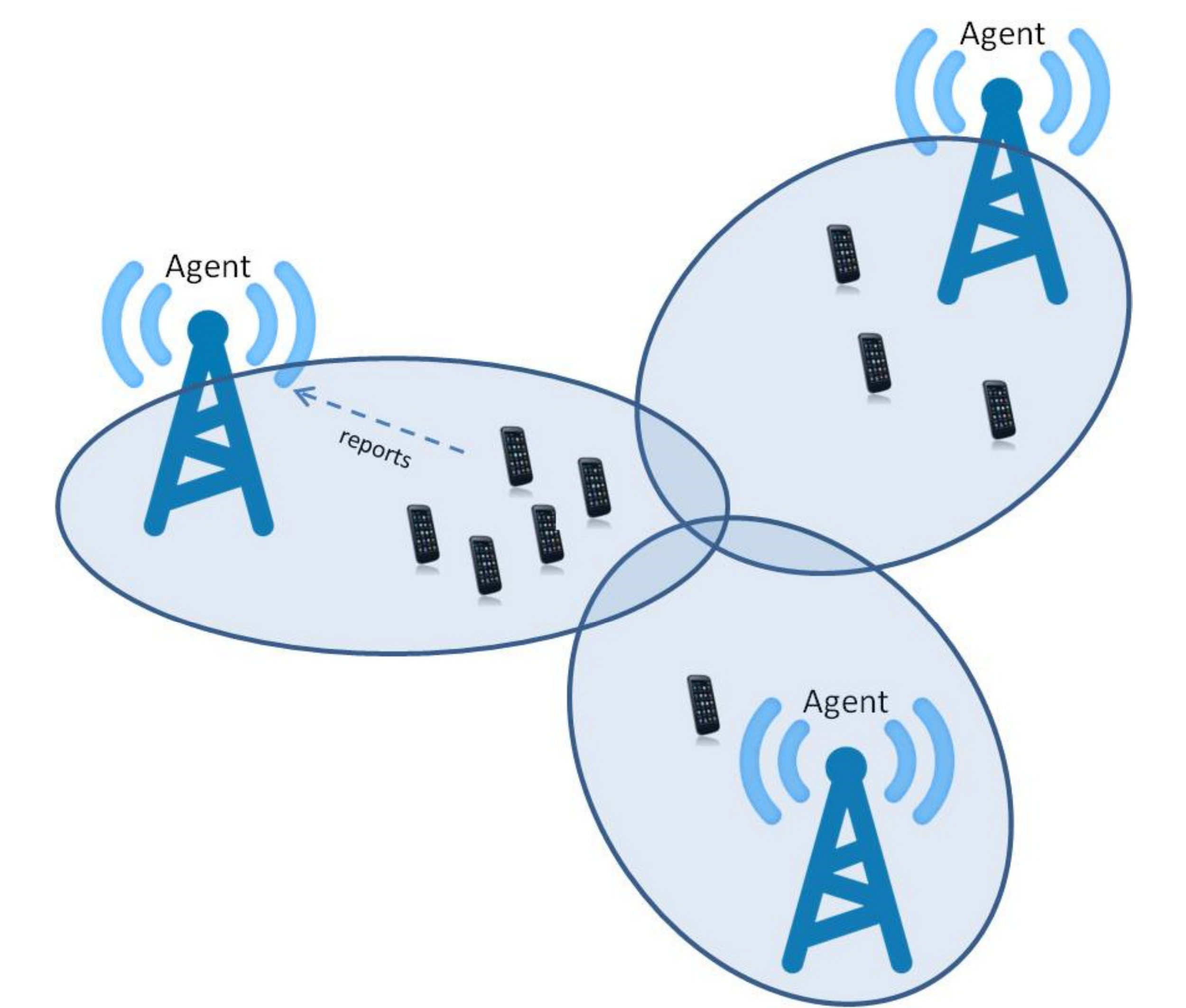}
}
\caption{An example of a cellular network with $3$ agents controlling the downlink transmission power of their own cell. }
\label{fig:RL_arch}
\end{figure}

\subsection{Reinforcement Learning}
Reinforcement Learning RL is an area of machine learning concerned with
learning how a software agent should learn to behave in a given environment in order to
maximize some form of cumulative reward.

The ingredients used to model the interaction of the agent with the
environment are the \emph{state}, the \emph{actions} and the \emph{reward}.
The state $s$ is a tuple of values, known as features, that describes the
agent in relation to the environment in a way that is relevant for the
problem at hand.
The action $a$ (chosen in the set of available actions) represents the
change (e.g., parameter changes) that the agent applies to
the environment in order to achieve his goal of maximizing the
given notion of reward.
The reward $r$ is a multi-objective scalar function which expresses what
we want the agent to achieve.

The agent therefore cycles through a series of transitions which consist of going from one state to the
next state by applying actions to the environment and receiving rewards as a consequence of his actions.
The evolution of the agent in the environment can therefore be described in terms of quadruples of the form
$(s_t, a_t, r_{t+1}, s_{t+1})$ where $t$ represents time.

The goal of RL is to derive a policy $\pi$, that is, a function that, given a state,
provides the action to take in order to maximize the cumulative reward:
\begin{equation}\label{eqn:policy_eqn}
\begin{aligned}
\pi(s_t) = a_t.
\end{aligned}
\end{equation}

Unlike Supervised Learning, where input data and output labels are given and a
function has to be derived from it, in RL the inputs are the state-action pairs, while the
outputs are the rewards (or, more precisely, a function of the reward) that the agents receives over time.

The rewards are often few and far in between and the relevant
actions that lead to certain outcomes might have been taken far back in
time and/or space. A central problem of RL is thus
that of properly assigning the reward to the actions that lead to such
reward (a notion known as \emph{credit-assignment}).

Another main issue of learning behavioral policies in domains that are unknown at first (especially in critical systems like a cellular network) is that of efficiently bringing the agent
from a tabula-rasa state to a condition where the agent is acting as close to optimality as possible.
This notion is also known as \emph{regret minimization} and is closely related to the topic of trading off
\emph{exploration} of the environment (to sample previously unseen parts of the state-action space) with
\emph{exploitation} of the knowledge accumulated so far.

A tabula-rasa agent that enters the world would therefore need to
explore the environment by applying actions and to gather data which are
as informative as possible so that a policy, encoding the knowledge
accumulated so far, is derived via a credit assignment algorithm. Over
time the agent will gradually move from a situation where it is mostly
exploring the environment to a situation where it is mostly exploiting
the accumulated knowledge. This gradual transition from an exploration strategy
to an exploitation strategy is referred to as epsilon-greedy. Every time an agent
has to take an action, it will take a random action with probability $\epsilon$ and a policy action with probability $1-\epsilon$.
To value of $\epsilon$ is gradually reduced from an initial value $\epsilon_{\max}$ to a final value $\epsilon_{\min}$ over
a predefined number of actions. For a complete treatment of the subject of RL, we refer the interested reader to \cite{SuB:98}.

\subsection{Q-learning}
Different classes of RL algorithms existing in literature are often
classified as critic-only, actor-only and actor-critic.

One of the most popular algorithm in the RL literature is that of Q-learning \cite{Wat:89},
which is a critic-only algorithm. Q-learning aims at learning an
action-value function (the so called Q-function) which gives the
expected utility of taking an action $a$ in a given state $s$ and following the policy $\pi$ afterwards.
Once such a function is learned, the policy can be derived from it by
evaluating the Q-function for every action and choosing the action
which gives the largest Q-value.

The Q-function for a given policy $\pi$ is defined as:
\begin{equation}\label{eqn:Qvalue_function}
\begin{aligned}
Q^{\pi}(s, a) = \mathbb{E}\left[\sum_{t=0}^\infty \gamma^t r(s_t, \pi(s_t))\vert s_0 = s, a_0 = a\right],
\end{aligned}
\end{equation}
The term $\gamma\in [0,1)$ is a discounting factor causing the value of rewards to
decay exponentially over time indicating the preference for immediate rewards while
making the optimization horizon for our agent finite (that is, the sum of rewards is finite).
In other words the Q-function learns to predict the expected cumulative discounted reward from taking
action $a$ in state $s$ and following the policy $\pi$ afterwards.

As the agent explores the environment by applying actions and receiving rewards,
it collects and stores transitions $(s_t, a_t, r_{t+1}, s_{t+1})$  in a growing batch.
From the batch of transitions a training sequence of input-output pairs is formed and
used to learn the Q-function.
The input is represented by state-action pairs $(s, a)$, while the output is represented
by the Q-function expressed as 

\begin{equation}\label{eqn:Qvalue_function}
\begin{aligned}
Q(s_t, a_t) = r_{t+1} + \gamma  \max_{a'}Q(s_{t+1}, a')
\end{aligned}
\end{equation}

Despite the fact that the target Q-function we are trying to learn contains
the Q-function itself, the presence of the reward term 
is sufficient for the target to be improved over time and move closer to the "real" Q-function.

This approach was originally conceived in conjunction with the usage of tables to store
the transitional data, but it becomes quickly unusable as the state-action space grows
and the required number of data samples becomes prohibitively large.
A more powerful approach consists then in combining Q-learning with function approximators which,
besides accelerating the learning, offer the possibility to generalize it to previously unseen states.

In our work we chose to represent the Q-function via an Artificial Neural Network (ANN).
The interested reader can find a detailed description of ANNs, their characteristics and
how they are trained in \cite{Bet:96}.

In the context of this paper suffices to say that we have used an ANN with 3 hidden
layers trained using the R-PROP algorithm \cite{RiH:93}
to update the weights of the ANN because it is consistent with our approach of using
the full growing batch of data to update the Q-function.

In short, to estimate the gradient we minimize the mean-square error
\begin{equation}\label{eqn:training_loss}
L=\dfrac{1}{2}\sum_t \vert r_{t+1}+\underset{a^\prime}{\max}\, Q(s_{t+1},a^\prime )-Q(s_t,a_t)\vert^2
\end{equation}
over the full batch of samples and then update the weights using the gradient descent with time-varying adaptive step sizes.

\subsection{Specifics of RL in the Context of Cellular Networks}

In the context of this paper, RL can be mapped to the entities in the network as follows:
\begin{itemize}
\item The environment is represented by the cellular network;
\item The agents are the logical network nodes implementing the RL algorithms and therefore capable of extracting
  a control policy (in our case power control) from the collected set of transitions;
\item The state, in its entirety, could be represented by the type and number of terminals, their traffic,
  their positions, their capabilities, the type and number of cells, the different measurements and KPIs, etc.
  but in general a more synthetic representation is used, that is, a more limited number of elements is added
  to describe the state of the agent in relation to the network;
\item The actions are typically represented in form of parameter adjustments (e.g., power-up, power-down, power-hold);
\item The reward could be represented as a function (non necessarily linear) of different Key Performance Indicators (KPI)s (capacity, coverage, delay, etc.).
\end{itemize}
The cellular network is usually affected by a few complications which are not present
in the basic formulation of RL.
\begin{enumerate}
\item The state as seen by the agent is only partially observable, that is, the agent can only observe a limited part of the state;
\item Once the action $a_t$ has been applied to the environment, the transition to the next state $s_{t+1}$ is
stochastic and not deterministic, that is $\mbox{Pr}(s_{t+1}=s'\vert s_t, a_t)$ for all $s'$;
\item In the environment there are several agents whose actions interfere with each other.
\end{enumerate}

Each of these problems is more or less severe depending on how much the actions of one agent influences the environment
as perceived by the other agents.
If the impact of the actions of each agent in the environment is such that the
reward signal is strongly corrupted then the learning will not be possible.
If this is not the case, the convergence rate will still be affected (because other
agents will still introduce noise into the system) but the agents will nonetheless be able to learn a policy.

As we show in the results of the paper, the latter is the case for downlink power control investigated here.
Each agent has access only to local information to construct the state but the reward is exchanged
across nodes and a network-wide reward function is constructed by aggregation of
individual reward values from the neighboring agents so as to encourage cooperation.
Moreover, in order to minimize the impact of other agents into the system, they coordinate
with each other by taking turns rather than acting all at the same time.

In the following, a more detailed description of each of the components of the RL framework
(that is, the state, action and reward) is provided.

\subsection{State}

In downlink power control, the state is represented by a set of features
constructed based on local measurements in each cell. These
features should provide concise information regarding the
performance of local cells as well as good indicators about
the situation of interfering cells. In particular, we employed the
following features to describe the state:
\begin{enumerate}
\item Cell power: One basic indicator of the state of the
network is indeed the amount of cell power. Since each
agent is placed in the base station, it has direct access into the
values of transmit power of each controlling cell.
\item Average Reference Signal Received Power (RSRP) in the cell: this feature is calculated by averaging
the RSRP measurements of users in the cell. The feature is essential since it implicitly embeds
an overview of users location in the cell. Cell-edge users
generally report a weak RSRP value. Thus, a low average RSRP in the cell is an indicator
of having majority of users at cell edge. 
\item Average interference in the cell: each user reports the received
signal strength from dominant interfering cells. The
average interference received from each of interfering
cells is then computed by taking the average of values
reported by users in the cell. Together with the average RSRP,
an agent can correlate the performance of 
its cell to the average channel quality of its served users
which in turn, is affected by the power values of other cells (and thus to the actions of other agents).
\item Cell reward: finally, the value of the reward is calculated in the
cell. It is used as part of the state so that its correlation to the global reward
can be used to provide a useful information for the
agent. The agent can then relate its own goal (local
reward) to the welfare of society (global reward).
\end{enumerate}

\subsection{Actions}

The true action space of downlink power control, the amount of cell power, is
a continuous variable. To avoid the action space explosion, the action space is quantized.
To strike a proper trade-off
between the speed of learning and the quality of the derived control
strategies, the action space should be reduced to a few discrete
values. In downlink power control, actions are such that agents can gradually
or more rapidly change their cells power level. A typical choice
in our simulations was given by the set
$\{0, \pm1, \pm3\}$dBs.

\subsection{Reward}

Since the objective of the RL problem is to maximize the cumulative discounted reward,
the design of the reward function is the only tool we have to enforce a particular agent
behavior compared to another.
Given that the wireless network is a cooperative environment (cells are not competing with
each other) it makes sense to define the reward in a global network-wide sense as function
of some fundamental KPIs.

For this work, we have considered the following general definition of $\alpha$-fair resource
allocation utility function \cite{AAG:08}.
\begin{equation}\label{eqn:reward_function}
\begin{aligned}
r(x) = \dfrac{1}{1-\alpha} \sum_{x_i\in \mathcal{X}} w_i (h_i(x_i)^{1-\alpha}-1), \quad \alpha \in [0,\infty),
\end{aligned}
\end{equation}
where $\alpha=[0,\infty)$ is a scalar coefficient,
and $h_i:\mathcal{X}\rightarrow \mathbb{R}$ is a transfer function. 
Here,  $x_i$ represents a radio measurement or a performance indicator associated with
the radio cell, $\mathcal{X}$ is the set of all radio measurements or performance
indicators associated with the radio cell and used for the definition of the
performance measurement, $w_i$ is a weight associated with $x_i$,
and $x=[x_1, \dots, x_{\lvert \mathcal{X}\rvert} ]$ is a vector comprising all $x_i\in \mathcal{X}$.

When the function $h_i (x_i )$ represent the average data throughput of user $i$ in the cell, for instance, the reward $r_t (x)$ in \eqref{eqn:reward_function} can be approximated for different values of $\alpha$ and weights $w_i$ with the following expressions:
\begin{enumerate}
	\item The average  throughput associated to the user devices in the cell, i.e., $r (x)=1/{\vert \mathcal{X}\vert} \sum_{x_i\in \mathcal{X}} h_i(x_i)$, if $\alpha=0$, $w_i=1/{\vert \mathcal{X} \vert}$ $\forall \,i$.
	\item The average  throughput associated to the cell, i.e., $r (x)= \sum_{x_i\in \mathcal{X}} h_i(x_i)$, if $\alpha=0$, $w_i=1$ $\forall \,i$.
	\item The average log-throughput of the users in the cell, i.e., $r (x)={ \sum_{x_i\in \mathcal{X}} \log(h_i(x_i))}/{\vert \mathcal{X}\vert}$, if $\alpha=1$, $w_i=1/{\vert \mathcal{X} \vert}$  $\forall \,i$.
	\item The average sum of log-throughput of the cell, i.e., $ r (x)=\sum_{x_i\in \mathcal{X}} \log(h_i(x_i))$, if $\alpha=1, w_i=1$ $\forall \,i$.
	\item The harmonic-mean throughput of users in the cell, i.e., \small$r (x)={\vert \mathcal{X}\vert}/ \sum_{x_i\in \mathcal{X}} h_i(x_i)^{-1}$, if $\alpha=2$, $w_i=1/{\vert \mathcal{X} \vert}$  $\forall \,i$.
	\item The harmonic-mean  throughput of the cell, i.e., $r (x)=1/{\sum_{x_i\in \mathcal{X}}} {h_i(x_i)}^{-1}$, if $\alpha=2$, $w_i=1$  $\forall \,i$.
\end{enumerate}
Each reward expression enables the agent node to optimize a different performance metric that can either be associated to individual user devices, to radio cells, or the whole cellular radio network.

\section{Performance evaluation}\label{sec:numerical}

We evaluate our reinforcement learning approach for power control and rate adaptation using a LTE-A compliant event-driven system-level simulator. Table \ref{tab:simulation_setup} summarizes the main simulation parameters.

\subsection{A 2-agent validating example}
We first evaluate the convergence our RL framework on a $2$-cell/2-agents example with $10$ UEs with full-buffer traffic deployed randomly in the network and compare its performance against the optimal one.
Agents take turns every $100$ milliseconds to control the donwlink power budget of the corresponding cell.
In this example, we model the cell state using a limited set of features comprising the cell's power, the average RSRP of the UEs in the cell, and the average interference measured at the UEs in the cell. We also evaluate the RL approach using two reward functions: the  harmonic-mean throughput of the network and the sum-log throughput of the network. The RSRP, interference, and reward functions are averaged with measurements collected over the period between two consecutive actions of an agent (i.e., 200ms). Starting from an exploration probability $\epsilon_{\max}=0.9$, each agent gradually annihilates this value until the minimum $\epsilon_{\min}=0.1$.


\begin{figure}[h!]
\centering{
\includegraphics[width=.9\linewidth]{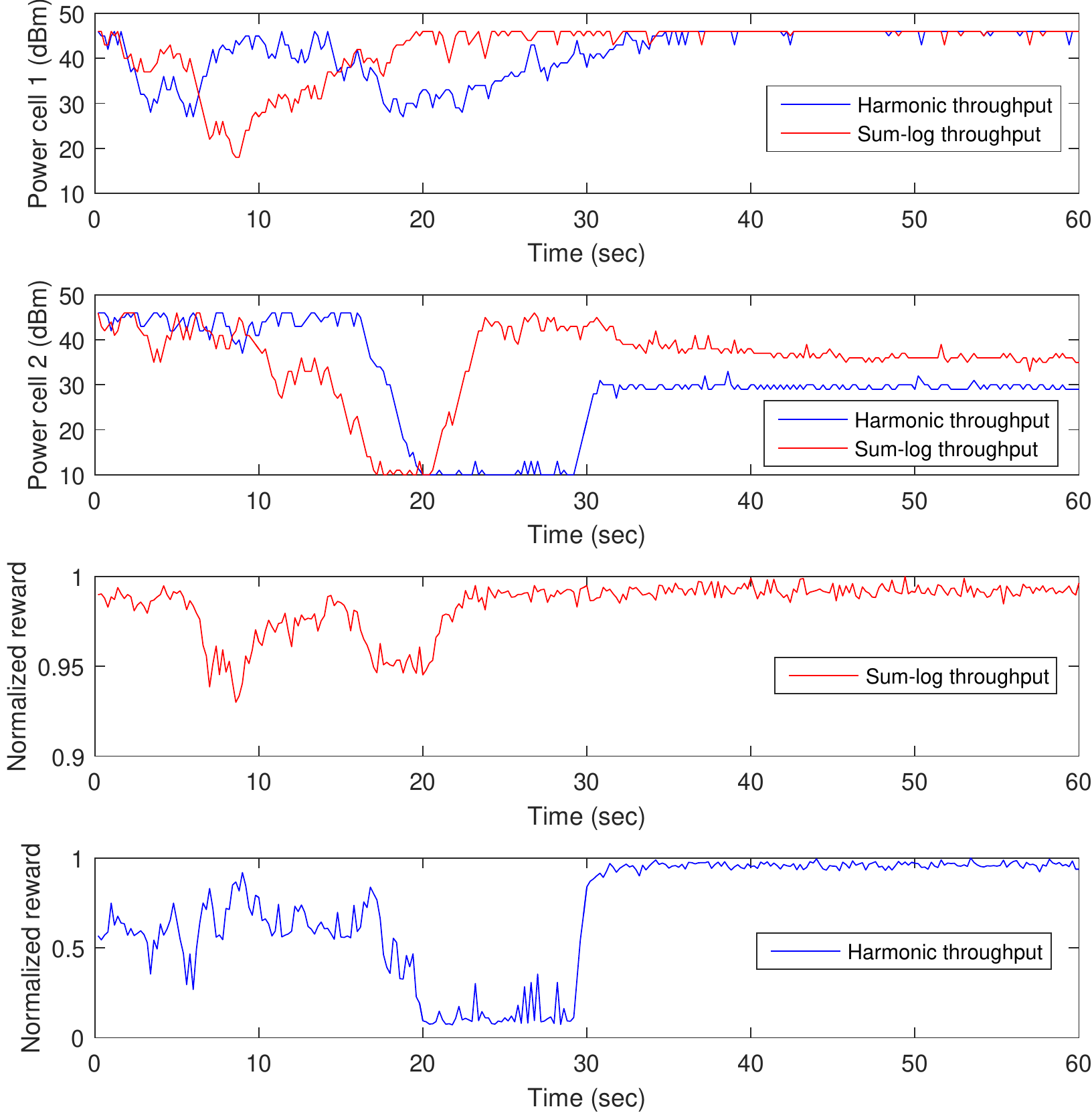}
}
\caption{ DL power level adjustments of two agents and associated normalized reward functions.}
\label{fig:validation_convergence}
\end{figure}

\begin{figure}[t!]
\centering{
\begin{subfigure}[t]{0.49\linewidth}
\centering
\includegraphics[width=\linewidth]{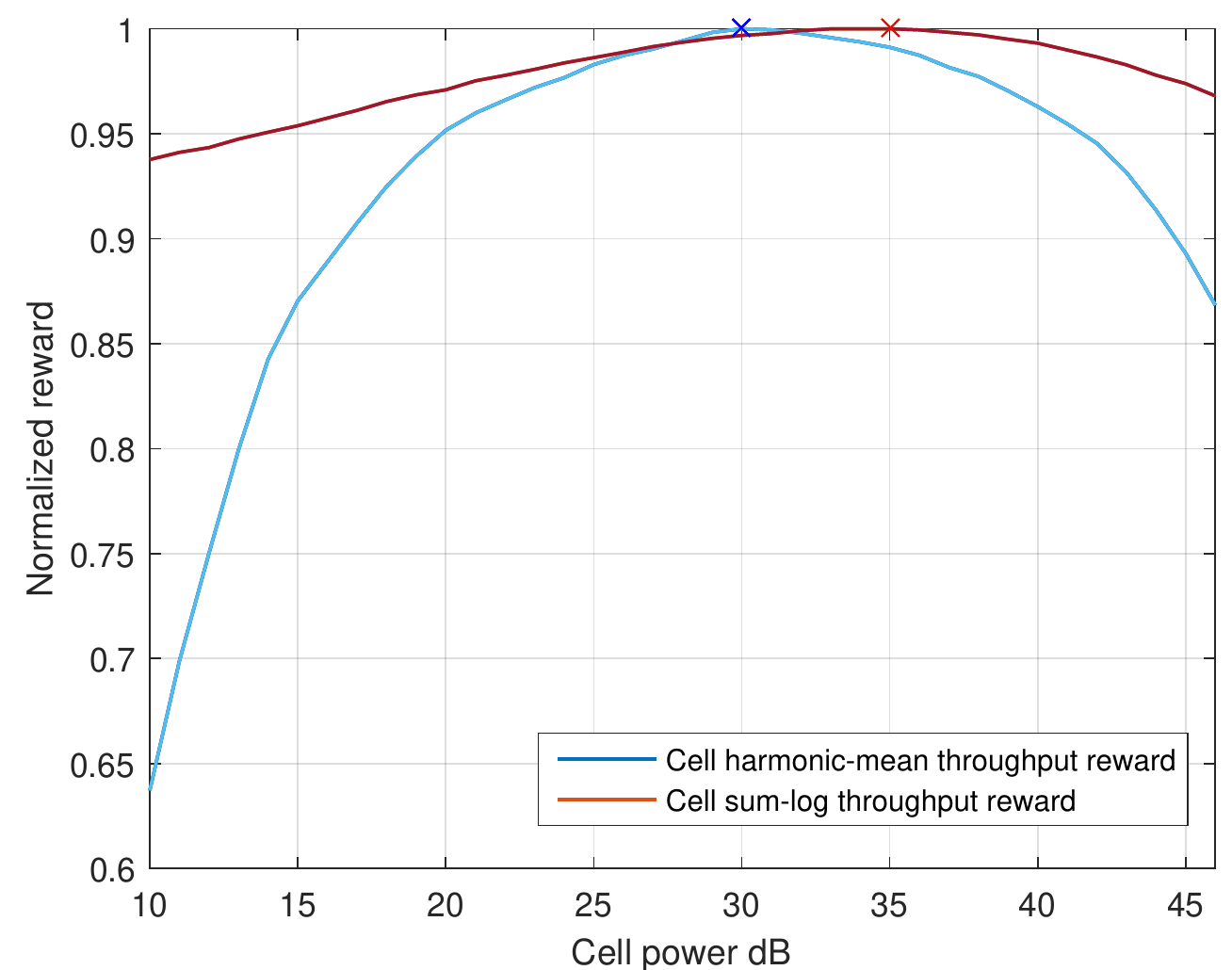}
\caption{Normalized reward functions.}
\label{fig:validation_HThroughput}
\end{subfigure}
~
\begin{subfigure}[t]{0.47\linewidth}
\centering
\includegraphics[width=\linewidth]{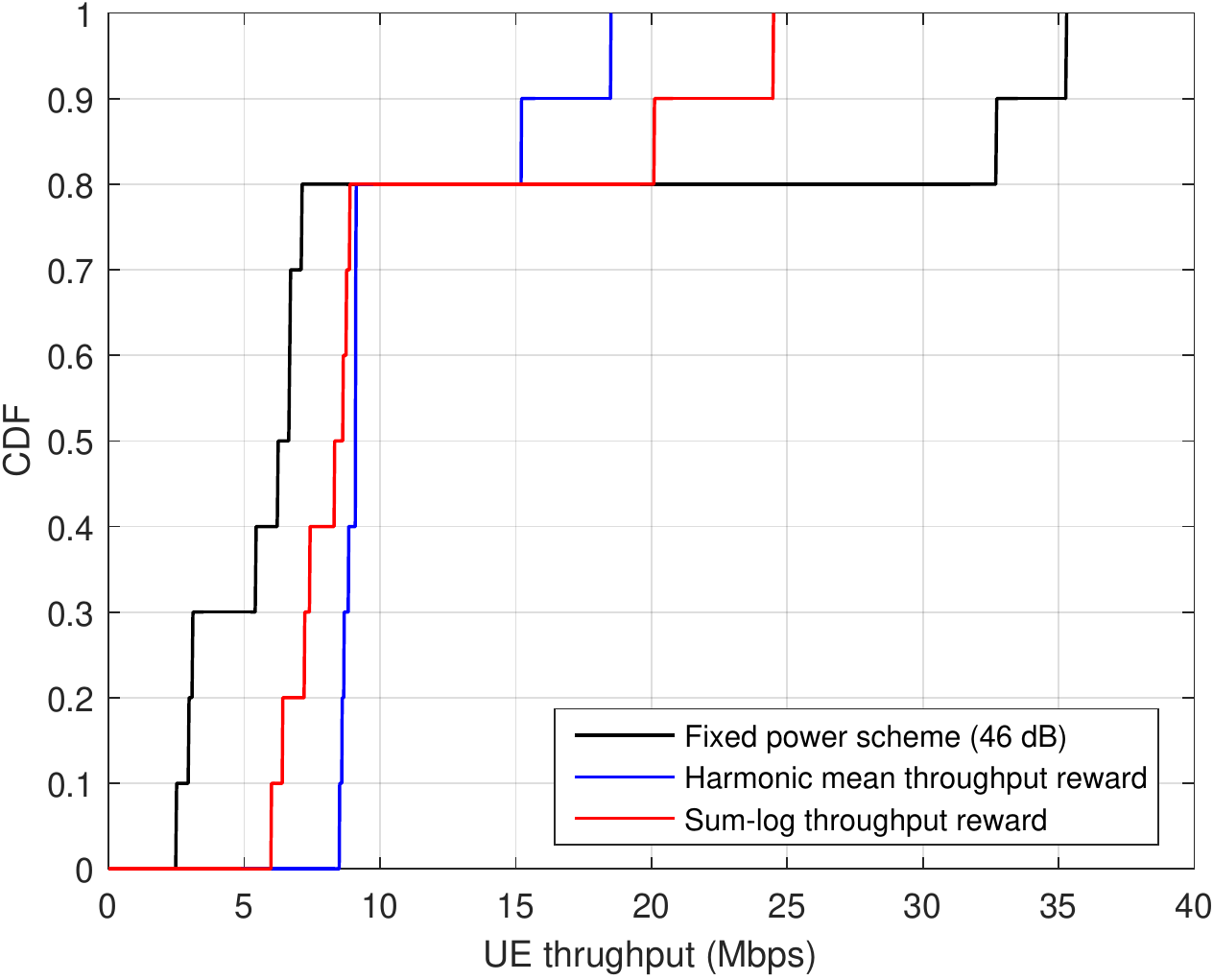}
\caption{CDF of UEs data rate.}
\label{fig:validation_cdf_thru}
\end{subfigure}
}
\caption{The left hand side shows the normalized reward functions for the two-cell example with $P_1 = 46$ dBm (optimal value) and $P_2$ ranging in $[10, 46]$ dBm with step size 1dBm. The right hand side shows the CDF of the user rates for the power values $(P_1^\star, P_2^\star)$ optimizing the two reward functions.}
\end{figure}

Fig.~\ref{fig:validation_convergence} shows the result for a 60s simulation. After an initial phase of intensive exploration, both agents quickly converge to a power control policy that maximize the chosen aggregate network reward.
%
%
%
We have then validated these results against the power levels optimizing both reward functions by exhaustively running a set of independent simulations with fixed power levels at both cells ranging in $[10,\; 46]$dBm with step-size 1dB\footnote{With a static system simulator or assuming average (fixed) channel gains over the simulation, the optimizer can be found solving problem~\eqref{eqn:NUMtransf}}.
%
Fig~\ref{fig:validation_HThroughput} shows that the RL algorithm converges to a power control policy optimizing (in average) the corresponding reward function\footnote{For ease of representation of both reward functions, Fig~\ref{fig:validation_HThroughput} shows results wherein the downlink power of cell 1 is fixed to 46dBm (the optimal value)}.
Fig~\ref{fig:validation_cdf_thru} compares the cumulative distribution function (CDF) of the user rates (averaged over the last 5s of simulation) for the case of fixed transmission power with 46dBm at both cells and for the power values achieved by the RL approach with the two reward functions. The results indicate that our design can converge to the optimal  power control policy that brings energy gain and fairness among users in the system, where the design of the reward function enables to tradeoff different degrees of fairness. 


\subsection{Extended simulations setup}
Since the gain from a dynamic downlink power allocation is assumed to be an effect of load balancing between the cells, we consider a 3-cell scenario with uneven load. In particular, three cells with relative load of $10\%$, $20\%$ and $70\%$, and a total number of $900$ UEs generated over 30 simulation drops.
We use the first half of the drops for training and the second half for evaluation\footnote{Training is also performed in the drops used for policy evaluation}.
%
We adopt the harmonic mean network throughput as reward as it offers fairness among the users and it strongly relates to the packet delay measure.  A summary of the results is given in Table \ref{tab:simulation_summary}, which shows the gain of the RL approach over an equal power allocation of $46$ dBm.

\begin{table}[t!]
\renewcommand{\arraystretch}{1.3}
\caption{Simulations parameters}
\centering
\begin{tabular}{lc}
    \toprule
    \textbf{Parameter}  &  \textbf{Value}\\
    \bottomrule
    \toprule
    TTI   &   $1$ms\\

    Bandwidth   &   $10$ MHz\\

    Default transmit power at BS   &   $46$ dBm\\


    Traffic type   &   Full buffer or bursty\\

    \#Cells (agents)   &   $2$ or $3$\\

	Reward  functions &   Harmonic-mean or sum-log throughput\\

	Action period  &  100 ms\\

	Agents scheduling  &   Round-robin\\

	Policy update period  &  $50$ data samples\\

	$(\epsilon_{\min}, \epsilon_{\max})$ & $(0.1, 0.9)$\\

	$\gamma$  &  $0.7$\\



\bottomrule
\end{tabular}
\label{tab:simulation_setup}
\end{table}

\begin{table}[t!]
\renewcommand{\arraystretch}{1.3}
\caption{Summary of simulation results for RL approach.}
\centering
\begin{tabular}{c c c c}
    \toprule
    \textbf{Scenario}  &  \textbf{$5\%$ (UE rate) gain} & \textbf{Median gain} & \textbf{Power reduction}\\
    \bottomrule
    \toprule
	Full buffer& $63\%$ & $15\%$ & $86\%$\\
Burst& $94\%$ & $22\%$ & $91\%$\\
    \bottomrule
\end{tabular}
\label{tab:simulation_summary}
\end{table}

\subsection{Full buffer traffic}
\begin{figure}[h]
\centering{
\includegraphics[width=0.6\linewidth]{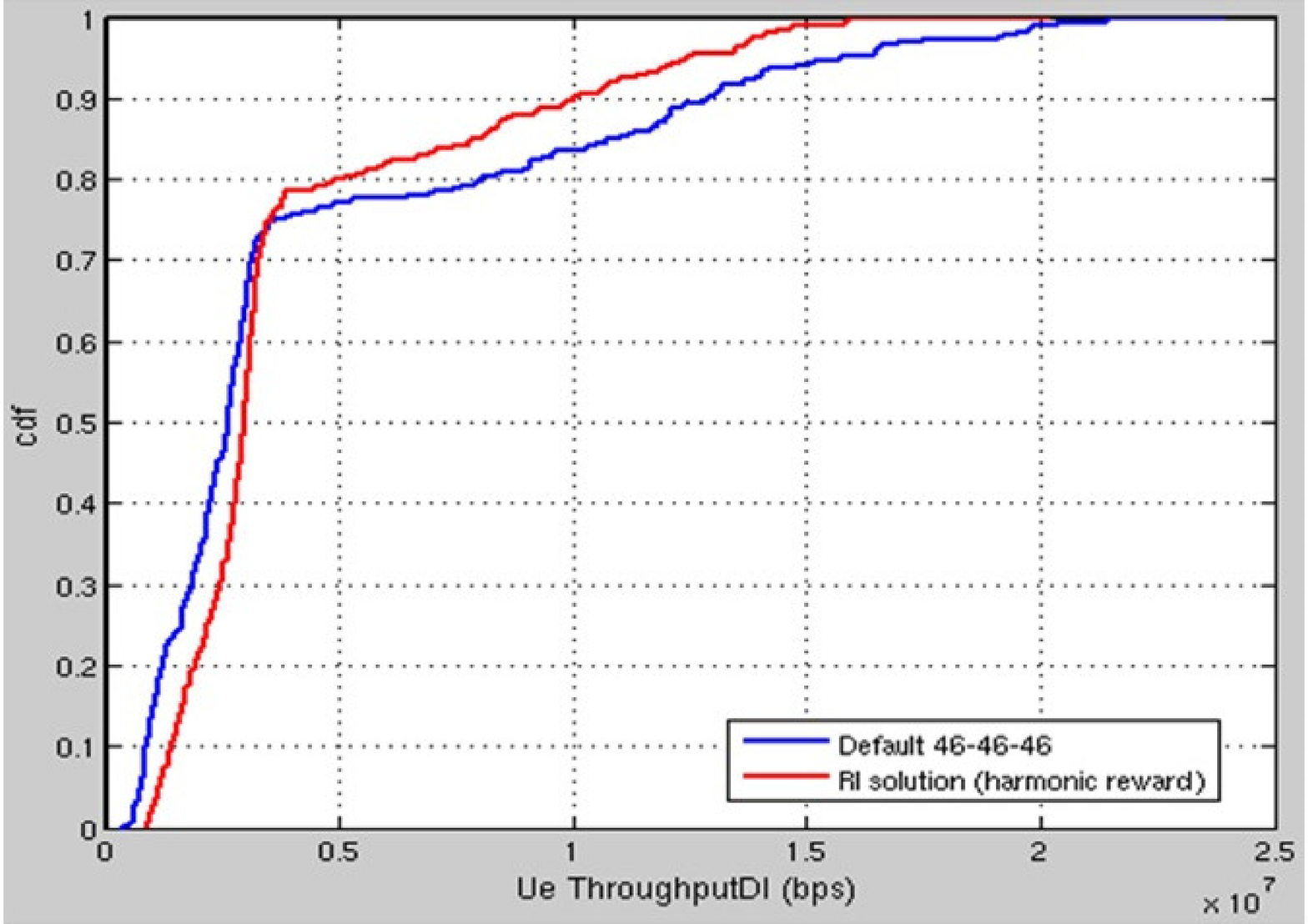}
}
\caption{CDF of the average user throughput with full buffer traffic.}
\label{fig:cdf_fullbuffer}
\end{figure}

\begin{figure}[h]
\centering{
\includegraphics[width=0.49\linewidth]{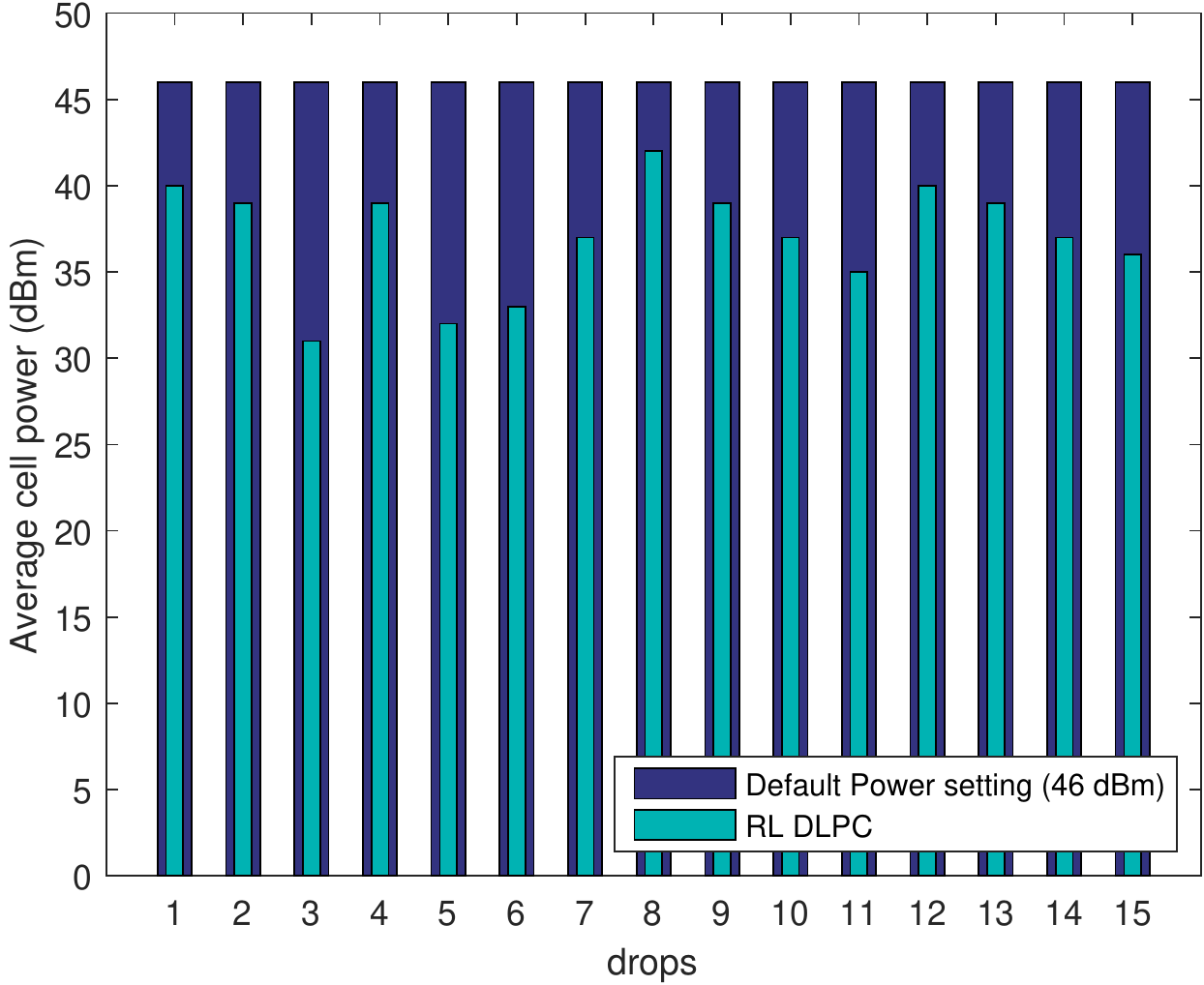}
\includegraphics[width=0.49 \linewidth]{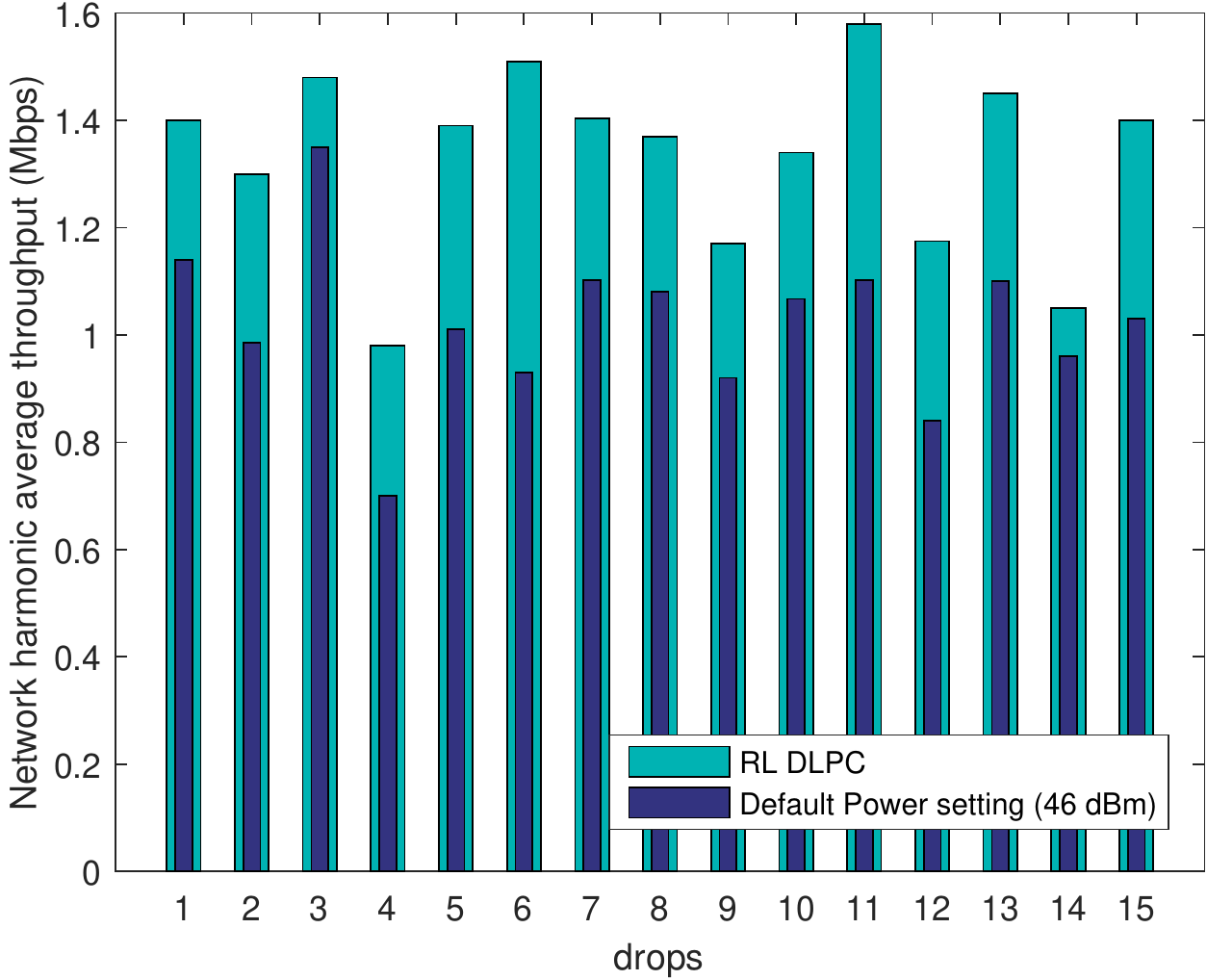}
}
\caption{Average cell transmit power and harmonic-mean network throughput with full buffer traffic.}
\label{fig:power_reward_fullbuffer}
\end{figure}

For full buffer traffic, we have simulated the training and evaluation drops for $60$s and $20$s, respectively, corresponding to $4000$ data samples per agent. Figure~\ref{fig:cdf_fullbuffer}, Figure~\ref{fig:power_reward_fullbuffer} and Table~\ref{tab:simulation_setup} summarize the results: compared to a baseline with fixed power allocation, the RL approach driven by an harmonic-mean reward brings fairness and throughput gains to low-rate users with a significant transmission power reduction.

Figure~\ref{fig:cdf_fullbuffer} show that the RL algorithm enforces inter-cell user fairness by reducing data rate to users in low loaded cells through a reduction of the associated downlink transmission power. In particular, the throughput gain ranges from $63\%$ for $5\%$-tile cell edge UEs (and up to the $20\%$-tile UEs) to a $15\%$ gain for the median UEs, respectively. Fairness, with full buffer traffic, comes at the price of a slight decrease in network throughput - below $10\%$ over all evaluation drops.
Figure~\ref{fig:power_reward_fullbuffer} show that RL approach further achieves an average downlink power reduction of $8.54$ dB corresponding to over $86\%$ power savings compared to the baseline, while effectively maximizing the reward function in each evaluation drop.

\subsection{Bursty traffic}
\begin{figure}[h]
\centering{
\includegraphics[width=0.6\linewidth]{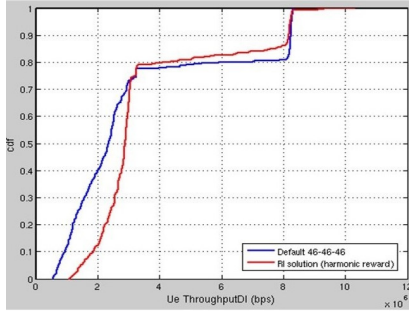}
}
\caption{CDF of users DL throughput for evaluation scenarios including $450$ UEs with bursty traffic.}
\label{fig:cdf_bursty}
\end{figure}

\begin{figure}[h]
\centering{
\includegraphics[width=0.49\linewidth]{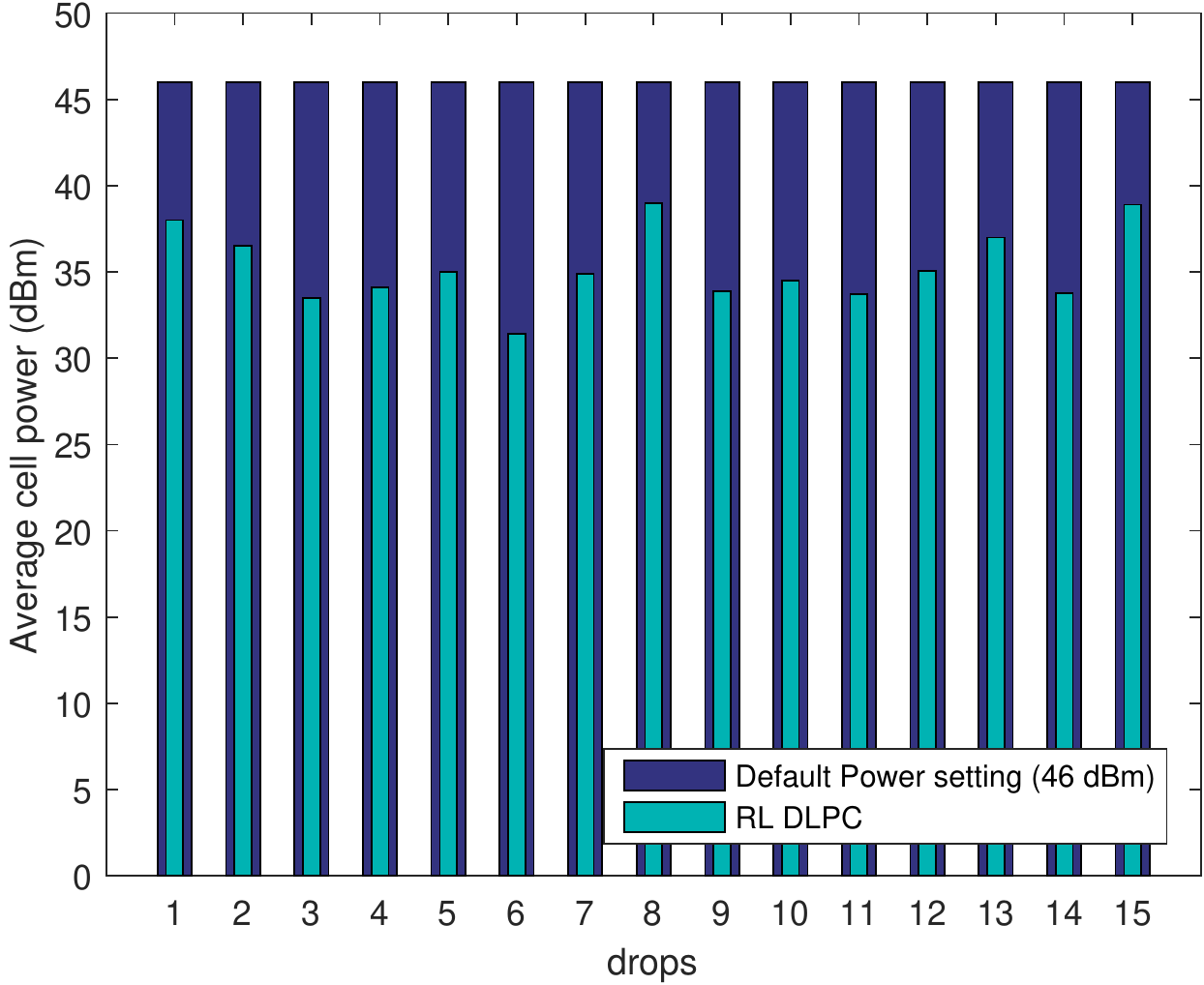}
\includegraphics[width=0.49 \linewidth]{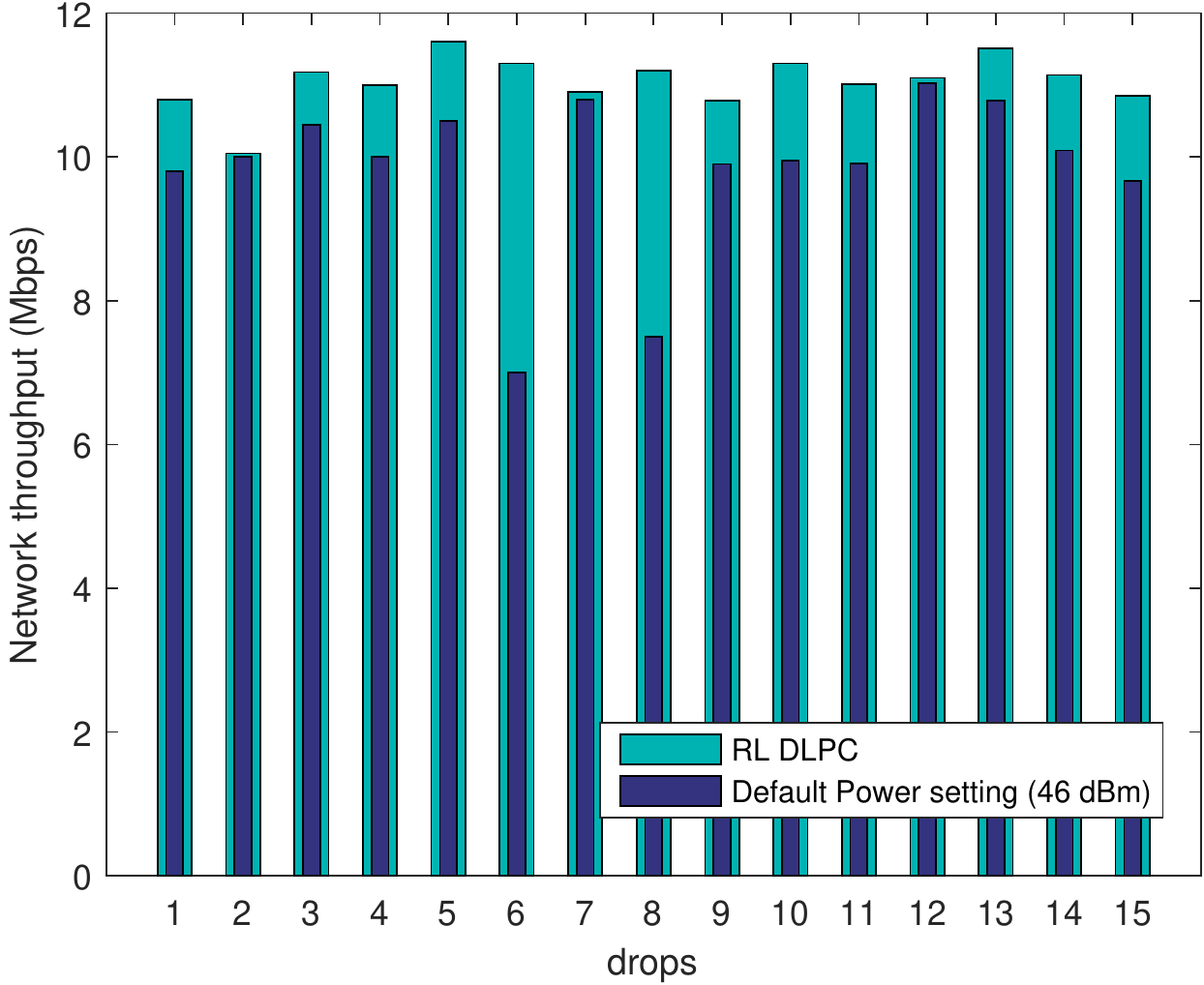}
}
\caption{Average cell transmit power and average network throughput for $15$ evaluation scenarios with bursty user-traffics.}
\label{fig:power_thru_bursty}
\end{figure}

For bursty traffic, we consider UEs downloading files of size $0.1$ MB with mean reading time $100$ms, corresponding to an average traffic  $1$ MB/s per UE. We have simulated the training and evaluation drops for $180$s and $10$s, respectively, corresponding to $9000$ training data samples per agent.
Figure~\ref{fig:cdf_bursty}, Figure~\ref{fig:power_thru_bursty} and Table~\ref{tab:simulation_setup} summarize the results: with bursty traffic, the RL approach successfully enforces inter-cell user fairness through a significant reduction of the cells transmit power without degrading the overall network throughput compared to a baseline with fixed transmit power.

Figure \ref{fig:cdf_bursty} shows that the RL approach achieves $94\%$ throughput gain for the $5\%$-tile UEs and  $22\%$ gain for median UEs. Overall, RL approach converges to a power control strategy that enables nearly $70\%$ of UEs  to enjoy a higher data rate.
Figure~\ref{fig:power_thru_bursty} shows an average transmit power reduction of $10.57$dB among the evaluation drops, i.e. over $91\%$ power saving, as well as an average network throughout gain of $14\%$ compared to the baseline.

\section{Conclusions}\label{sec:conclusion}
In this paper, we considered the problem of downlink power control for cellular systems. In particular, we provided a RL based method that adapts the power budget of cells to the dynamic conditions of the network and user traffics.  With proper features and reward functions selection, we presented simulation results where tuning cell powers using our algorithm offers significant improvements over baseline for both full and bursty traffic scenarios.
An interesting future direction is to further investigate the multi-agent aspects of  learning framework by increasing the number of cells and making the power control action user specific rather than cell specific.
\balance

\bibliography{ICCreferences}

\end{document}